\documentclass{amsart}

\usepackage{amsfonts}
\usepackage{verbatim}
\usepackage{amssymb}
\usepackage{amscd}
\usepackage{amsmath}
\usepackage{amsthm}
\usepackage{graphicx}
\usepackage{subfigure}
\usepackage{geometry}
\usepackage{ifpdf}
\usepackage[parfill]{parskip}
\ifpdf
\fi
\makeatletter
\def\thm@space@setup{%
  \thm@preskip=\parskip \thm@postskip=0pt
}
\makeatother  

\makeatletter
\def\subsubsection{\@startsection{subsubsection}{3}%
  \z@{.5\linespacing\@plus.7\linespacing}{.1\linespacing}%
  {\normalfont\itshape}}
\makeatother

\makeatletter
\def\subsection{\@startsection{subsection}{3}%
  \z@{.5\linespacing\@plus.7\linespacing}{.1\linespacing}%
  {\normalfont\itshape}}
\makeatother

\newtheorem{theorem}{Theorem}[section]

\newtheorem{lemma}[theorem]{Lemma}

\newtheorem{corollary}[theorem]{Corollary}

\begin{document}
\bibliographystyle{plain}

\title{A generalization of the Birthday problem}
\author{Sukhada Fadnavis}
\email{Sukhada Fadnavis <sukhada@math.harvard.edu>}
\address{Dept. of Mathematics, Harvard University, One Oxford Street, Cambridge, MA 02138.}
\date{\today}
\maketitle

\begin{abstract} The birthday paradox states that there is at least a 50\% chance that
some two out of twenty-three randomly chosen people will share the same birth date. The calculation for this problem assumes that all birth
dates are equally likely. We consider the following two modifications of this question. If the distribution of birthdays is non-uniform, does that increase or decrease the probability of matching birth dates?
Further, what if we focus on birthdays shared by some particular pairs rather than any two people. Does a non-uniform distribution on birth dates increase or decrease the probability of a matching pair? In this paper we present our results 
in this generalized setting. We use some results and methods due to Sokal \cite{SokalZeroes} concerning bounds on the roots of chromatic polynomials to prove our results. 
\end{abstract}

\section{Introduction}\label{Introduction}

The Birthday problem is a classical and well-studied problem in elementary probability. There is a vast literature on this problem and it's
generalizations and their applications; for example see \cite{VonMises}, \cite{Stein}, \cite{CamarriPitman}, \cite{PersiSusan},
\cite{DiaconisMosteller}, \cite{Holst}. The birthday problem asks for the minimum number $n$ of birthdays that we need to sample independently
so that the probability that all of them are distinct is small (say less than 50\%). The well known answer to this question is 23. To see this,
suppose we have $n$ people each having one of $q$ possible birthdays distributed uniformly and independently. The probability that everybody
has a distinct birthday is:
\begin{equation}
    \prod_{i=1}^{n-1}\left(1- \frac{i}{q}\right).
\end{equation}
For $q = 365$ this probability goes below 0.5 for the first time when $n=23$. \\

One wonders though if it is accurate to assume that all birthdays occur with equal probability. There are more induced births during the
weekdays than on weekends because of ready availability of staff. There may be fluctuations in birthrates during different seasons. Does this
affect the probability of two students sharing the same birthday? If so, does the probability increase or decrease? It is known (for example,
see \cite{Uniformbest1}, \cite{Uniformbest2}, \cite{Uniformbest3}) that the probability of matching birthdays increases if the distribution of
birthdays is not uniform. To see this, let $\textbf{p}= (p_1, \ldots, p_q)$ be the distribution on the $q$ birth dates and let $P_n(p_1,
\ldots, p_q)$ denote the probability that no two people share the same birthday under this distribution. Then,
\begin{equation}
    P_n(p_1, \ldots, p_q) = n! \sum_{i_1 < \ldots < i_n}(p_{i_1}\ldots p_{i_n}).
\end{equation}

By a classical theorem of Muirhead \cite{Muirhead} this is a concave symmetric function of the $p_i's$. Hence,
\begin{equation}\label{uniformbest}
    P_n(p_1, \ldots, p_q) \leq P_n\left(\frac{1}{q}, \ldots \frac{1}{q}\right).
\end{equation}
Thus, in this case the uniform distribution is the worst case distribution i.e. the probability of all distinct birth dates is maximixed when the birthdates are uniformly distributed.  \\

Further generalizing the situation, what happens if instead of all distinct birth dates we just want all pairs of friends to have distinct
birth dates? We construct a friendship graph $G$ as follows: there is a vertex corresponding to each person and an edge between two if and only
if they are friends. Now replacing birth dates by $q$ colors we get the following graph theory problem. Consider a graph $G$ on $n$ vertices.
Suppose the vertices are colored at random with $q$ colors occurring with probabilities $p_1 \cdots p_q$. We say that a coloring of a graph is
a proper coloring if no edge is monochromatic. Let $P_G(p_1, \ldots, p_q)$ denote the probability that the random coloring thus obtained is a
proper coloring. In this setting the Birthday Problem asks for the smallest $n$ such that,
\begin{equation}
    P_{K_n}\left(\frac{1}{q}, \ldots, \frac{1}{q}\right) \leq \frac{1}{2}.
\end{equation}

In the general setting the distribution $\textbf{p}= (p_1, \ldots, p_q)$ need not be uniform. Also $G$ can be any underlying graph which we
call the friendship graph. Equation \eqref{uniformbest} tells us that $P_{K_n}(p_1, \ldots, p_q)$ is maximized if all the colors occur with
probability $p_i = 1/n$, where $K_n$ denotes the complete graph on $n$ vertices. A natural question to ask is if this is true for all
underlying graphs $G$, i.e.
\begin{equation}\label{question}
    \text{Is } P_G(1/q, \ldots, 1/q) \geq P_G(p_1, \ldots, p_q) \text{ for all graphs } G?
\end{equation}
The answer to this question is negative as shown by the following example due to Geir Helleloid:

\textbf{Example (Geir Helleloid):} Consider the `star graph' $K_{1,4}$ colored with two colors $c_1, c_2$ with respective probabilities $p_1,
p_2$. Here $P(\frac{1}{2}, \frac{1}{2}) = \frac{1}{2^4}$. On the other hand $P(\frac{1}{5}, \frac{4}{5}) = \frac{4^4}{5^5} + \frac{4}{5^5} >
\frac{1}{2^4}$. In general if $G = K_{1,n}$ for $n\geq 4$, then,
\begin{equation}
    P_G\left( \frac{1}{2}, \frac{1}{2} \right) < P_G\left( \frac{1}{n+1}, \frac{n}{n+1} \right).
\end{equation}

Note that as we increase $q$ the situation changes. In fact we will show  in Section \ref{Star graphs} that for star graphs $G = K_{1,n}$ the
probability $P_G$ is indeed maximized by the uniform distribution when $q \geq n$.

\begin{figure}[htp]
\centering
\includegraphics[width=2in]{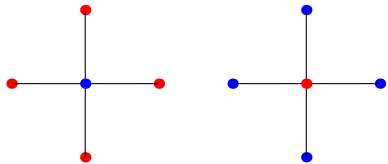}
\caption{Four star and it's two proper colorings with two colors.} \label{K14}
\end{figure}

\begin{figure}[htp]
\centering
\includegraphics[width = 2.5 in]{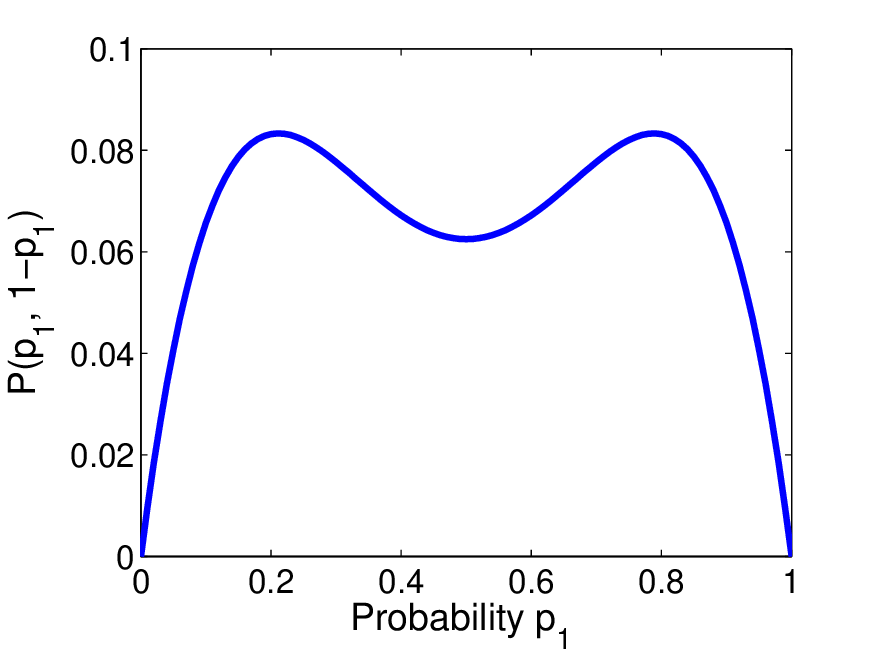}
\label{star} \caption{Above is a plot of $P_{K_{1,4}}(p_1, 1-p_1)$ against $p_1$. We see that $P_{K_{1,4}}(p_1, 1-p_1)$ is maximized at 1/5 and
4/5. }
\end{figure}

In this paper we show that such counterexamples can exist only for `small' values of $q$. If $q$ is large in terms of the maximum degree of the graph, then the answer to question  \ref{question} is positive. More precisely, we have the following theorem: 

\begin{theorem}\label{MainTheorem} If $G = (V,E)$ is a graph with maximum degree $\Delta$, then for $q> 4 \times 10^4 \Delta^4$ we have,
\begin{equation}
    P_G\left(\frac{1}{q}, \ldots, \frac{1}{q}\right) \geq P_G(p_1, \ldots, p_q),
\end{equation}
for any distribution $\textbf{p} = (p_1, \ldots , p_q)$ on the colors.
\end{theorem}

The following special cases were studied in \cite{Fadnavis}:

\begin{theorem}[\cite{Fadnavis}] \label{clawfreeproof} If $G$ is claw-free then $P_G(p_1, \ldots p_q)$ is maximized when $p_1 = \cdots = p_q = 1/q$. In fact $P_G$ is Schur-concave on the set of probability distributions $\textbf{p} = (p_1,\ldots, p_q)$.
\end{theorem}

\begin{theorem}[\cite{Fadnavis}] \label{shameful} If $G = (V,E)$ is a graph with maximum degree $\Delta$, then for $q> 36\Delta^{3/2}$ we have,
\begin{equation}
    P\left(\frac{1}{q-1}, \ldots, \frac{1}{q-1}\right) \leq P\left(\frac{1}{q}, \ldots, \frac{1}{q}\right)
\end{equation}
\end{theorem}

The remaining paper is organized as follows: In section \ref{Star graphs} we prove a stronger result for the special case of $G$ being a star graph. The proof of Theorem \ref{MainTheorem} is provided in the section \ref{Proof2}. The proof uses some results and methods due to Sokal \cite{SokalZeroes} concerning bounds on the roots of chromatic polynomials. \\

Before we conclude the introduction, we would like to point out how the chromatic polynomial is related to this problem:

\subsection{Graph coloring and chromatic polynomials}

Throughout this paper we will assume that $G = (V,E)$ is a finite simple graph on $n$ vertices with maximum degree $\Delta$. We say that a function $\alpha: V \rightarrow \{1, \ldots, q\}$ is a $q$-coloring of
$G$ if for each edge $(u,v)$ of $G$ we have $\alpha(u) \neq \alpha(v)$.  Let $\chi_G(q)$ be the number of $q$-colorings of $G$. In general
given a graph $G$ it is difficult to say whether it has a $q$-coloring or not, and hence it is also difficult to count the exact number of
$q$-colorings of $G$. Using inclusion exclusion we see that $P_G$ is in fact a polynomial known as the \emph{chromatic polynomial}:
\begin{equation}
    \chi_G(q) = \sum_{E' \subset E} (q)^{C(E')} (-1)^{|E'|},
\end{equation}
where $C(E') $ denotes the number of connected components in $E'$. \\

We note that $P_G(p_1, \ldots, p_q)$ can also be written as a polynomial of $p_1, \ldots, p_q$ in a similar manner:
\begin{equation}
    P_G(p_1, \ldots, p_q) = \sum_{E' \subseteq E} (-1)^{|E'|} \prod_{\substack{\gamma \subset E' \\ \gamma \text { connected.}}}  (p_1^{|\gamma|} + \ldots + p_q^{|\gamma|}),
\end{equation}
where the sum goes over all subsets $E' $ of the edge set $E$, and the product is over all connected components of $(V, E')$. By $|\gamma|$
we denote the number of vertices in $\gamma$. Note that the two polynomials are related to each other by the following equality:
\begin{equation}
    P_G\left(\frac{1}{q}, \ldots, \frac{1}{q} \right) = \frac{\chi_G(q)}{q^{n}}.
\end{equation}

Due to this similarity the study of $P_G(p_1, \ldots, p_q)$ is similar to the study of the chromatic polynomial $\chi_G(q)$. This is useful because the chromatic polynomial is a very well-studied object. The literature on chromatic polynomials is vast and we refer the reader to \cite{Read}, \cite{DongBook} for excellent surveys. For the purposes of this paper we will be interested in the study of the roots of the chromatic polynomial \cite{Brown}, \cite{BRW}, \cite{SokalZeroes}, \cite{Borgs}. 

\section{Star graphs}\label{Star graphs}
Before we proceed with the proof in the general case, let us first consider the case of the star graph. In this case we have the following result: 
\begin{theorem}\label{Stars} For the star graph $G = K_{1,n}$ and $q>n$ we have,
\begin{equation}
     P_G(p_1, \ldots, p_q) \leq P_G\left( \frac{1}{q}, \ldots, \frac{1}{q} \right).
\end{equation}
\end{theorem}

\begin{proof}
 Given the star graph and colors as above, the probability that a random coloring gives rise to a proper coloring is:

\begin{equation}
    P_G(p_1, \ldots, p_q) = \sum_{i=1}^{q} p_i(1-p_i)^n.
\end{equation}

\begin{figure}[htp]
\centering
\includegraphics[width = 3 in]{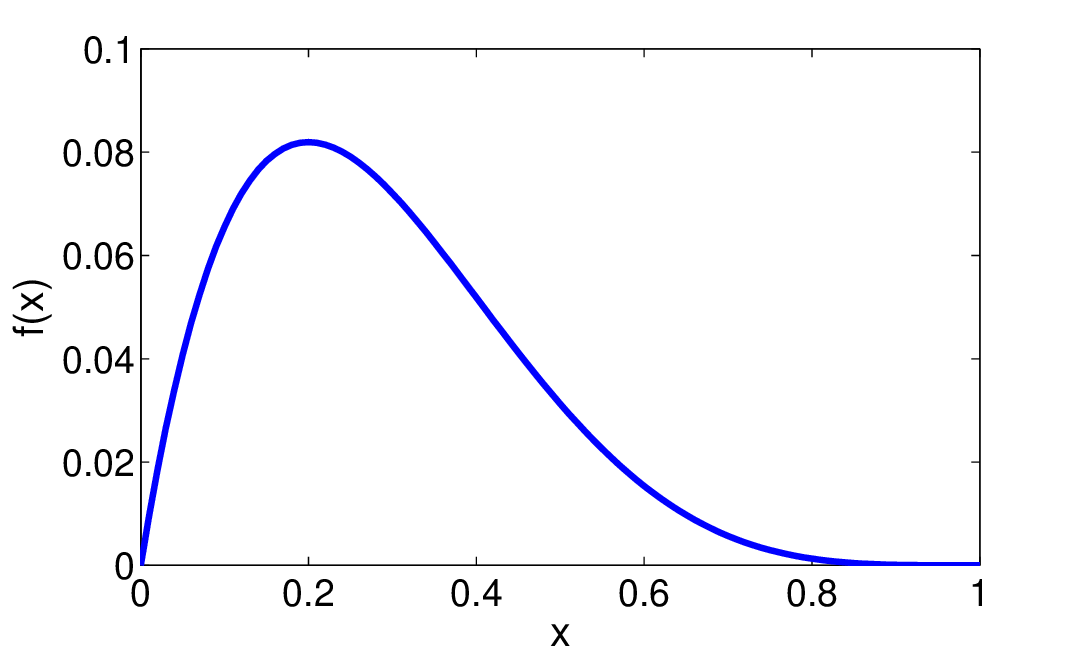}
\caption{This is a plot of $f(x) = x(1-x)^4$ against $x$. We see that $f$ is maximized at 0.2 and is concave on $[0,0.4]$. } \label{star}
\end{figure}

Note that the function $f(x) = x(1-x)^n$ is unimodal for $0 \leq x \leq 1$. In fact, it is concave on $[0,\frac{2}{n+1}]$ and convex on
$[\frac{2}{n+1},1]$. The function has a unique maxima at $\frac{1}{n+1}$ on the interval $[0,1]$. Let $$\Omega = \{x_1 \ldots x_q | x_i \geq 0,
x_1+ \ldots x_q = 1\}.$$ We wish to show that $P_G$ has a maximum at $(\frac{1}{q}, \ldots ,\frac{1}{q})$, on $\Omega$. Let $$\Theta = \{ (x_1,
\ldots, x_q) \in \Omega | x_i \leq \frac{2}{n+1} \text{ for all } i \}.$$ Then by the unimodality and concavity of $x(1-x)^n$ on
$[0,\frac{2}{n+1}]$, it follows that $P_G$ has a maxima at $(\frac{1}{q}, \ldots ,\frac{1}{q})$, on $\Theta$. Now suppose $(x_1, \ldots, x_q)
\in \Omega$ is such that $x_i > \frac{2}{n+1}$ for some $i$. Then there is also an $x_j$ such that $x_j < \frac{1}{n+1}$. Then replacing $x_i$
by $x_i + x_j - \frac{1}{n+1}$ and $x_j$ by $\frac{1}{n+1}$ increases the value of $P_G$. Continuing thus, we can get to a point in $\Theta$
where the value of $P_G$ will be strictly greater than the value of $P_G$ at the point outside $\Theta$ where we started. This together with
the earlier fact proves that $P_G$ has a maximum at $(\frac{1}{q}, \ldots ,\frac{1}{q})$, on $\Omega$.
\end{proof}

\section{Proof for general graphs (Proof of Theorem \ref{MainTheorem})}\label{Proof2}

\begin{proof}
As in the case of the star graph, the proof in the general case has two steps. The first step is to show that if any $p_i$ is much larger than
$1/q$ then, $P_G(1/q, \ldots, 1/q) \geq P_G(p_1, \ldots, p_q)$. More precisely,

\begin{theorem}\label{Step1}(Proved in \ref{step1proof}) If $\displaystyle{p_i \geq 2\sqrt{\frac{\Delta}{q}}}$ for some $i$, then $P(p_1, \ldots, p_q) \leq P(1/q, \ldots, 1/q)$.
\end{theorem}

The next step is to show that when all the $p_i's$ are close to $1/q$ then $P_G$ is log-concave for large enough $q$:

\begin{theorem}\label{Step2}(Proved in \ref{step2proof}) If $q > 4\times 10^4 \Delta^4$, then $P_G(p_1, \ldots, p_q)$ is maximized at  $(1/q, \ldots, 1/q)$ in the region  
$$ \Omega = \left\{(p_1, \ldots, p_q)  \in \mathbb{R}_{+}^q: p_{1}^m + \ldots + p_{q}^m \leq \left( 2\sqrt{\frac{\Delta}{q}}\right)^{m-1} \forall m \in \mathbb{Z}_{+} \right\} $$ 
\end{theorem}

We note one small lemma before completing the proof the theorem. 
\begin{lemma}\label{discbound} Let,
    \begin{equation}
        \Omega_1 = \left\{(p_1, \ldots, p_q) : p_i \geq 0,  p_1 + \ldots + p_q =1,  p_i \leq 2\sqrt{\frac{\Delta}{q}}  \right\}
    \end{equation}
and,
\begin{equation}
   \Omega = \left\{(p_1, \ldots, p_q)  \in \mathbb{R}_{+}^q: p_{1}^m + \ldots + p_{q}^m \leq \left( 2\sqrt{\frac{\Delta}{q}}\right)^{m-1} \forall m \in \mathbb{Z}_{+} \right\},
\end{equation}
as above. Then, $\Omega_1 \subset \Omega$.
\end{lemma}
\begin{proof}
    Let, \begin{equation}
            \left \lfloor \sqrt{ \frac{q}{4\Delta}} \right \rfloor = k \text{ and } a = 1- 2k\sqrt{\frac{\Delta}{q}}  \leq 2\sqrt{\frac{\Delta}{q}}.
            \end{equation}
            Since $\Omega_1$ is a symmetric convex polytope and $p_1^m + \ldots + p_q^m$ is a symmetric convex function it is maximized on the endpoints.
    Thus, $p_1^m + \ldots + p_q^m \leq k \left( 2\sqrt{\frac{\Delta}{q}}\right)^m + a^m \leq \left( 2\sqrt{\frac{\Delta}{q}}\right)^{m-1}$
    since $a^m \leq ab^{m-1}$ for all $b \geq a \geq 0$.
\end{proof}

Theorem \ref{Step1}, Theorem \ref{Step2} and Lemma \ref{discbound} together prove Theorem \ref{MainTheorem}. We prove Theorems  \ref{Step1} and \ref{Step2} in the following sections. 
\end{proof}

\subsection{Proof of Theorem \ref{Step1}}\label{step1proof}

\begin{proof}
Let $N = \chi_G(q)$ be the number of proper colorings of $G$ using $q$ colors. Suppose the vertices of $G$ have degrees $d_1, \ldots, d_n$
respectively. Then $2|E| = \displaystyle\sum_{i = 1}^{n}d_i$. Note that for $q>\Delta$,
\begin{equation}\label{coloringineq}
\frac{N}{q^n} \geq \prod_{i \leq n}\left(\frac{q-d_i}{q}\right) \geq \displaystyle\left(\frac{q - \Delta}{q}\right)^{(\sum d_i)/\Delta} =
\left(1-\frac{\Delta}{q}\right)^{2|E|/\Delta}.
\end{equation}

The first inequality follows by coloring vertices in a fixed order. Vertex $i$ can have any of $q-n_i \geq q - d_i$ colors, where $n_i$ is the
number neighbors of vertex $i$ that have already been colored. To see the second inequality, note that for $1 \geq a \geq b \geq 0$ and
$\epsilon \geq 0$ one has,
\begin{equation}
    (1 - a -\epsilon)(1-b + \epsilon) = 1 -a -b + ab -\epsilon (a-b) - \epsilon^2 \leq (1-a)(1-b).
\end{equation}
This implies that $\log (\prod _{i \leq n}(1-x_i))$ is schur-concave. Thus,
\begin{equation}
    \prod_{i \leq n}\left(1- \frac{d_i}{q}\right)^{\Delta} \geq  \left(1-\frac{\Delta}{q}\right)^{2|E|} \times 1^{n \Delta - 2|E|},
\end{equation}
since $(d_1, \ldots, d_1, \ldots, d_n ,\ldots, d_n) \preceq (\Delta, \ldots, \Delta, 0, \ldots, 0)$ where the first vector has $\Delta$
co-ordinates that are $d_i$ for each $i$ and the second vector has $2|E|$ co-ordinates that are $\Delta$ and the rest are 0's. This gives the
second inequality in \ref{coloringineq}.

Hence,
\begin{equation}
    P(1/q,\ldots ,1/q)=\frac{N}{q^n} \geq \left(1-\frac{\Delta}{q}\right)^{2|E|/\Delta}.
\end{equation}

Now since the maximum degree is $\Delta$ we can find a set $U \subset E$ of $\lceil |E|/2\Delta \rceil$ disjoint edges in $G$. Hence,
\begin{equation}
    P(p_1,\ldots, p_q) \leq (1-\sum p_i^2)^{|E|/2\Delta}
\end{equation}

So now it suffices to prove that
\begin{equation}
    (1-\sum p_i^2)^{|E|/2\Delta}  \leq \left(1-\frac{\Delta}{q}\right)^{2|E|/\Delta},
\end{equation}

that is,
\begin{equation}
    (1-\sum p_i^2) \leq \left(1 - \frac{\Delta}{q}\right)^{4}.
\end{equation}

Or, since

\begin{equation}
    1- \frac{4\Delta}{q} \leq \left(1 - \frac{\Delta}{q}\right)^{4},
\end{equation}

it suffices to prove that

\begin{equation}
    \begin{split}
        &(1-\sum p_i^2) \leq 1- \frac{4\Delta}{q} \\
        & \text{ i.e.\ }\frac{4\Delta}{q} \leq \sum p_i^2.
    \end{split}
\end{equation}

This is true by the hypothesis and hence completes the proof.
\end{proof}

\subsection{Proof of Theorem \ref{Step2}}\label{step2proof}

For the proof of Theorem \ref{Step2} we will make extensive use of ideas and theorems due to A. Sokal \cite{SokalZeroes} and C.Borgs
\cite{Borgs}. The first hurdle is to get a nice combinatorial, inductive formula for $P_G$. As stated earlier, inclusion-exclusion gives:
\begin{equation}
    P_G(p_1, \ldots ,p_q) = \sum_{E' \subseteq E} (-1)^{|E'|} \prod_{\gamma \in  \mathcal{C}(E')} (p_1^{|\gamma|} + \ldots + p_q^{|\gamma|}),
\end{equation}
where $\mathcal{C}(E')$ denotes the set of all connected components $\gamma$ of $(V, E')$ and by $|\gamma|$ we denote the number of vertices in $\gamma$. Also
note that the summand is 1 when $E' = \emptyset$. To see this, recall that if  $A = A_1 \cup \ldots \cup A_k$ is a union of events then inclusion exclusion
gives:
\begin{equation}
    \mathrm{Prob}(A) = \sum_{i \leq k} \mathrm{Prob}(A_i) - \sum_{1 \leq i < j \leq k} \mathrm{Prob}(A_i \cap A_j) + \ldots + (-1)^{k+1}\mathrm{Prob}(A_1 \cap \ldots \cap A_k).
\end{equation}
So, let $A$ be the event that the coloring is not a proper coloring and let $A_i$ denote the event that edge $i$ is monochromatic (i.e. both
end points have the same color). Then since $A = A_1 \cup \ldots \cup A_{|E|}$, and $P_G(p_1, \ldots, p_q) = 1- \mathrm{Prob}(A)$, we get,
\begin{equation}
    \begin{split}
            P_G(p_1, \ldots, p_q) &  =  1 - \sum_{\emptyset \neq E' \subseteq E} (-1)^{|E'| + 1} \prod_{\gamma \in \mathcal{C}(E')}  (p_1^{|\gamma|} + \ldots + p_q^{|\gamma|}) \\
        & = \sum_{E' \subseteq E} (-1)^{|E'|} \prod_{\gamma \in \mathcal{C}(E')}  (p_1^{|\gamma|} +\ldots+p_q^{|\gamma|}). \\
    \end{split}
\end{equation}

Thus, we can think of $P_G$ as a complex multivariate polynomial $P_G(z_1, \ldots, z_q)$. Now $P_G$ can be rewritten by collecting together
subsets $E'$ of $E$ that lead to connected components on the same set of vertices. Let $\mathcal{G} = (\mathcal{V}, \mathcal{E})$ denote the
graph whose set of vertices is given by the set of connected subsets $S$ of $V$ such that $|S| \geq 2$. There is an edge between $S_1$ and
$S_2$ if $S_1 \cap S_2 \neq \emptyset.$ Then, $P_G$ can be rewritten as:

\begin{equation}\label{inclusion-exclusion}
    \begin{split}
         & P_G(z_1, \ldots , z_q) =  \sum_{\substack{\mathcal{W} \subseteq \mathcal{V} \\ \mathcal{W} \text{ independent}}} \prod_{S_i \in S} w(S_i) \\
         & \text{where } w(S) = (z_1^{|S|} + \ldots + z_q^{|S|}) \sum_{\substack{\gamma \subseteq E,\\ (S,\gamma) \text{ connected}}} (-1)^{|\gamma|},
    \end{split}
\end{equation}
where the summand is $1$ when $\mathcal{W} = \emptyset$.

One advantage of writing $P$ in this form is that it can be decomposed nicely. Let $\mathcal{U} \subseteq \mathcal{V}$. We define:
\begin{equation}
    P_{\mathcal{U}} = \sum_{\substack{\mathcal{W} \subseteq \mathcal{U} \\ \mathcal{W} \text{ independent}}} \prod_{S_i \in S} w(S_i).
\end{equation}

Let $\eta \in \mathcal{V}$, and let $\mathcal{V}' = \mathcal{V} \setminus \{\eta\}$. Further let, $\mathcal{V}_0 = \mathcal{V} \setminus N[\{\eta\}]$, where
$N[x]$ denotes the set containing $x$ and it's neighbors in $\mathcal{G}$.  Then,
\begin{equation}
    P_{\mathcal{V}} = P_{\mathcal{V}'} + w(\eta) P_{\mathcal{V}_0}.
\end{equation}

Such a decomposition is useful for proving statements inductively. For example, it is used to prove Dobrushin's theorem which gives conditions
under which functions, which can be decomposed as above, are non-zero. Applying a version of Dobrushin's theorem (as explained in section
\ref{logboundproof}) gives us the following result:

\begin{theorem}[Proved in \ref{logboundproof}]\label{log-bound} Let $\Delta$ be the maximum degree of $G$ and let $K = 7. 963907$ be a constant. If $q > K^2\Delta^3$ then $|\log P_G(z_1, \ldots, z_q)| \leq 4|E|/5$ in the region 
$$\Omega_0 = \left\{(z_1, \ldots, z_q)  \in \mathbb{C}^q: |z_1^m + \ldots + z_q^m| \leq \left( 2\sqrt{\frac{\Delta}{q}}\right)^{m-1} \forall m \in \mathbb{Z}_{+} \right\}.$$

\end{theorem}

The above theorem tells us that the Taylor expansion of $\log P_G(z_1, \ldots, z_q)$ converges in the region $\Omega_0$. The next theorem provides bounds on the the coefficients of this Taylor expansion. 

\begin{theorem}[Proved in \ref{logboundproof}]\label{CoefBounds}

In the above setup $\log P_G(p_1, \ldots, p_q)$ can be expressed as the power series of $\nu_{\mathbf{z}}(m)'$s where $$\nu_{\mathbf{z}}(m) = z_1^m + \ldots + z_q^m.$$ The expansion has the form,
\begin{equation}
    \log P_G(z_1,\ldots, z_q) = -|E|(z_1^2 + \ldots + z_q^2) + \sum_{M=3}^{\infty}\sum_{\substack{\alpha = (\alpha_1 \leq \ldots, \leq \alpha_s): \\ \sum \alpha_i = M, \  \alpha_i \geq 2}} C_{\alpha } \prod_{1}^{s} \nu_{\mathbf{z}}(\alpha_i),
\end{equation}
where $C_\alpha$ are constants. The series converges in $\Omega_0$. The first couple of coefficients are given by,

\begin{equation}\label{coefs}
C_{(2,2)} = -\sum_{i} \binom{d_i}{2}  \text{ and } C_{(3)} \leq \sum_{i} \binom{d_i}{2}.
\end{equation}
The remaining coefficients in the expansion are bounded above as follows,
\begin{equation}
    |C_{\alpha}| \leq \frac{4}{5}|E|\times \left(\frac{K\Delta}{2}\right)^{M-s}, \text{ when } \alpha_1+\ldots +\alpha_s = M.
\end{equation}

\end{theorem}

Finally, we need a small lemma before we complete the proof of Theorem \ref{Step2}

\begin{lemma}\label{replacement} Let $\Theta = \{(a_1, \ldots, a_q): \sum a_i = 1\}$. Let $f$ be a function on $\Theta$. If $$g(a_1, \ldots, a_q) = f(a_1, \ldots, a_q) - (a_1^{s+r} +\ldots + a_q^{s+r})$$
is minimized on $\Theta$ at $(1/q, \ldots, 1/q)$ then so is $$h(a_1, \ldots, a_q) = f(a_1, \ldots, a_q) - (a_1^{s+1} +\ldots +
a_q^{s+1})(a_1^{r}+ \ldots + a_q^{r}).$$ \end{lemma}
\begin{proof}
Note that, $$h(a_1, \ldots, a_q) = g(a_1, \ldots, a_q) - (a_1^{s+1} +\ldots + a_q^{s+1})(a_1^{r}+ \ldots + a_q^{r}) + (a_1^{s+r} +\ldots +
a_q^{s+r}).$$ Now, since $g$ is minimized at $(1/q, \ldots, 1/q)$, it suffices to prove that $$w(a_1, \ldots, a_q) = -(a_1^{s+1} +\ldots +
a_q^{s+1})(a_1^{r}+ \ldots + a_q^{r}) + (a_1^{s+r} +\ldots + a_q^{s+r})$$ is minimized at $(1/q, \ldots, 1/q)$. This is true since $w(1/q,
\ldots, 1/q) = 0$ and in general $w(a_1, \ldots, a_q) \geq 0$. To see this, note that,
\begin{equation}
     \begin{split}
        w(a_1, \ldots, a_q) =& -(a_1^{s+1} +\ldots + a_q^{s+1})(a_1^{r}+ \ldots + a_q^{r}) \\
        & + (a_1^{s+r} +\ldots + a_q^{s+r})(a_1+\ldots + a_q) \\
        & = \sum_{i \neq j}(a_i^{1}a_j^{s+r} + a_j^{1}a_i^{s+r} - a_i^{s+1}a_j^{r} - a_j^{s+1}a_i^{r}) \geq 0 \text{ by AM-GM }.
     \end{split}
\end{equation}
This completes the proof.

\end{proof}

Finally, in the proof of Theorem \ref{Step2} we use corollary \ref{CoefBounds} to show that when $q$ is large enough (as stated in the theorems) the first term of the Mayer expansion dominates  which further implies the result. 

\begin{proof}
As observed above,
\begin{equation}
    \begin{split}
        \log P_G(p_1,\ldots, p_q) &= -|E|(p_1^2 + \ldots + p_q^2) + C_{(3)} (p_1^3 + \ldots + p_q^3) -\sum_{i=1}^n \binom{d_i}{2} (p_1^2 + \ldots + p_q^2)^2 \\
        &+ \sum_{M=5}^{\infty}\sum_{\substack{\alpha = (\alpha_1 \leq\ldots \leq \alpha_s): \\ \text{ partition of } M \\ \alpha_i \geq 2}}C_{\alpha } \prod_{i\leq s} (p_1^{\alpha_i} + \ldots + p_q^{\alpha_i})
    \end{split}
\end{equation}

and,
\begin{equation}
	\begin{split}
    		& C_{\alpha} \leq \frac{4|E|}{5}\left(\frac{K\Delta}{2}\right)^M, \text{ for } \alpha \text{ a partition of } M, \\
		& C_{(3)} \leq \sum_{i=1}^n \binom{d_i}{2}.
	\end{split}
\end{equation}

Now, by Theorem \ref{replacement} it suffices to show that $\tilde{P}_G(p_1, \ldots, p_q)$ is maximized when $p_1 = \ldots = p_q$, where,

\begin{equation}
    \begin{split}
        \tilde{P}_G(p_1, \ldots, p_q) &=  -|E|(p_1^2 + \ldots + p_q^2) + C_{(3)} (p_1^3 + \ldots + p_q^3) \\
        &+ \sum_{M=5}^{\infty} \sum_{\substack{\alpha = (\alpha_1 \leq\ldots \leq \alpha_s): \\ \text{ partition of } M \\ \alpha_i \geq 2}}\frac{4|E|}{5}|V|\left(\frac{K\Delta}{2} \right)^M (p_1^{M-s+1} + \ldots + p_q^{M-s+1})\\
        & = C_{(3)} (p_1^3 + \ldots + p_q^3)  -|E|(p_1^2 + \ldots + p_q^2) \\ & + \sum_{k=3}^{\infty} A(k)\times \frac{4|E|}{5}\left(\frac{K\Delta}{2} \right)^{k}(p_1^{k+1} + \ldots + p_q^{k+1}),
    \end{split}
\end{equation}
where, $A(k)$ denotes the number of ordered partitions of $k$. The second equality follows since for every partition
$\alpha = (\alpha_1 \leq \ldots \leq \alpha_s)$ of $M$ such that $\alpha_i \geq 2$, we get a unique partition $\beta = (\alpha_1-1 \leq \ldots \leq \alpha_q-1) $ of $M-s$. Note, $A(k) \leq 2^k$. The Hessian of $\tilde{P}_G(p_1, \ldots, p_q)$ is a diagonal matrix with i'th
diagonal entry given by,
\begin{equation}
	\begin{split}
    		H_{ii} & = -2|E| + 6 C_{(3)} p_i + \frac{4|E|}{5}\sum_{k=3}^{\infty} A(k) \times \left(\frac{K\Delta}{2} \right)^{k}k(k+1)p_{i}^{k-1} \\
		& \leq -2|E| + 6 \sum_{i=1}^n \binom{d_i}{2} p_i + \frac{4|E|}{5}\sum_{k=3}^{\infty} A(k) \times \left(\frac{K\Delta}{2} \right)^{k}k(k+1)p_{i}^{k-1}
	\end{split}
\end{equation}
Since $\sum_i d_i = 2 |E|$ and $d_i \leq \Delta$, we have,
\begin{equation}
    \sum_{i} \binom{d_i}{2} \leq \frac{1}{2}\sum_i d_i^2 \leq |E|\Delta. 
\end{equation}
Using above inequality and $A(k) \leq 2^k$ gives,
\begin{equation}
    \begin{split}
        H_{ii} & \leq -2|E| + 6 \Delta |E| p_i + \frac{4|E|}{5}\sum_{k=3}^{\infty} 2^{k} k(k+1)\left(\frac{K\Delta}{2} \right)^{k}p_{i}^{k-1} \\
        & \leq -|E|\left(2 - 6\Delta p_i - \frac{4}{5} \sum_{k \geq 3} k(k+1)(K\Delta)^{k}p_{i}^{k-1}  \right) 
    \end{split}
\end{equation}
Using $p_i \leq \left( \frac{4\Delta}{q}\right)^{1/2}$ we get,
\begin{equation}
		H_{ii}  \leq -|E|\left( 1 - 6\Delta \left( \frac{4\Delta}{q}\right)^{1/2} + 1 - \frac{4}{5} \sum_{k \geq 3}  k(k+1)(K\Delta)^{k}\left( \frac{4\Delta}{q}\right)^{\frac{k-1}{2}}   \right) 
\end{equation}

Let $ x = K\Delta \left( \frac{4\Delta}{q}\right)^{1/2} < 1$. Then,

\begin{equation}
		H_{ii}  \leq -|E|\left( 1 - 6\Delta \left( \frac{4\Delta}{q}\right)^{1/2} + 1 - \frac{4}{5}K \Delta \frac{2x^2(3x^2-8x+6)}{(1-x)^3}  \right) 
\end{equation}

Recall that $7 < K <8$ and $\Delta \geq 1$. Thus, choosing $$x < \frac{0.1}{\Delta^{1/2}} \text{ that is, } q > 4 \times 10^4 \Delta^3 $$ gives that  $H_{ii} < 0.$ Further $H_{ii} < 0$ implies that $\tilde{P}$ is log-concave. Also, $\tilde{P}$ is symmetric in the $p_i$'s, hence log-concavity implies that it is minimized at $(1/q, \ldots, 1/q)$. This completes the proof of
Theorem \ref{MainTheorem}.

\end{proof}

\subsubsection{Proof of theorem \ref{log-bound}}\label{logboundproof}

In this section we will prove Theorem \ref{log-bound}. We will need the following theorems due to A. Sokal \cite{SokalZeroes} and C.Borgs \cite{Borgs}. First we explain some notation and then state three equivalent versions of Dobrushin's theorem, which we will use in the proof. \\

Let $X$ be a set (called a `single particle state space') with relation $\sim$ on $X\times X$ and and $w : X \rightarrow \mathbb{C}$  a complex function called the \emph{fugacity vector}. \\
We say $X' \subseteq X$ is \emph{independent} if $ x \sim y $ for all $x, y \in X'$.\\
Let,
\begin{equation}
        Z_X(w) = \sum_{\substack{X' \subseteq X \\ X' \text{ independent}}} \prod_{x \in X'}w_x.
\end{equation}

\begin{theorem}[Dobrushin's theorem as stated in \cite{Borgs}]
In the above setup $Z_X$ is non-zero in the region $|w_x| \leq R_x$, if  there exist constants $c_x \geq 0$ such that,
\begin{equation}
    R_x \leq (e^{c_x} - 1) \exp\left(-\sum_{y \nsim x}c_y \right).
\end{equation}
Further,
\begin{equation}
    \left| \log\left\{\frac{Z_X}{Z_{X'}}\right\} \right| \leq \sum_{x \in X \setminus X'} c_x, \text{ for all } X' \subseteq X.
\end{equation}
Hence, in particular,
\begin{equation}
    |\log Z_X| \leq \sum_{x \in X} c_x.
\end{equation}
\end{theorem}

From Dobrushin's theorem follows the Kotecky-Preiss condition: \\

\begin{theorem}[Kotecky-Preiss condition] In the above setup $Z_X$ is non-zero in the region $|w_x| \leq R_x$, if  there exist constants $c_x \geq 0$ such that,
\begin{equation}
    R_x \leq c_x \exp\left(-\sum_{y \nsim x}c_y \right)
\end{equation}
\end{theorem}

We will use the following consequence of the Kotecky-Preiss condition as stated by Sokal \cite{SokalZeroes},\\

\begin{theorem}[Proposition 3.2 of \cite{SokalZeroes}]\label{Dobrushin} Let $R_x \geq 0$ for all $x \in X$. Suppose that $X = \bigcup_ {n = 1}^{\infty} X_n$ is a disjoint union such that there exist constants $\{A_n\}_{n=1}^{\infty}$ and $\alpha$ such that,
\begin{enumerate}
\item  $\displaystyle \sum_{y\in X_n : y\nsim x} R_y \leq A_{n}m $ , for all $m,n$ and all $x \in X_m$.

\item $\displaystyle \sum_{n=1}^{\infty} e^{\alpha n}A_n \leq \alpha.$
\end{enumerate}

Then the Kotecky-Preiss condition holds with the choice $c_x = e^{\alpha n}R_x$ for all $x\in X_n$.\\

\end{theorem}

\begin{corollary}\label{alphabound} Assume the hypothesis of Theorem \ref{Dobrushin}. Further let $F \subseteq X_2$ be such that for all $y \in X_n$ there is a $v \in F$ such that $y\nsim v$. Then,
\begin{equation}
    |\log Z_X| \leq \sum_{x \in X} c_x \leq |F|\alpha.
\end{equation}
\end{corollary}
\begin{proof}
By choosing $m=2$ in part 1 of Theorem \ref{Dobrushin} we have,
\begin{equation}
        \sum_{y\in X_n : y\nsim v} e^{\alpha n}R_y \leq 2 e^{\alpha n} A_{n} \text{ for all } v \in F.
\end{equation}

Thus,
\begin{equation}
        \sum_{x \in X} c_x \leq \sum_{n \geq 1} \sum_{v \in F}  \sum_{\substack{y\in X_n \\ y\nsim v}} e^{\alpha n}R_y \leq \sum_{n \geq 1}\sum_{v \in F}2e^{\alpha n} A_{n} \leq \sum_{v \in F} \sum_{n \geq 1} 2e^{\alpha n} A_{n} \leq 2|F|\alpha.
\end{equation}
The last inequality follows from condition 2 of Theorem \ref{Dobrushin}.
\end{proof}

Next, we state four theorems that were used in \cite{SokalZeroes} to prove a bound on the roots of the chromatic polynomial. We will use these results in a very similar fashion in our proof. \\

\begin{theorem}[Penrose's Theorem \cite{Penrose}]\label{Penrose}
Let $G = (V,E)$ be a finite graph on $n$ vertices. Then,
\begin{equation}
\left|\sum_{\substack{E'\subseteq E \\ (V,E') \text{connected}}} (-1)^{|E'|}\right| \leq T_n(G),
\end{equation}
where $T_n(G)$ denotes the number of spanning trees of $G$.\\
\end{theorem}

\begin{theorem}[Special case of Proposition 4.2 in \cite{SokalZeroes}] \label{spanningbound} Let $H$ be a graph degree $\Delta$ and let $x$ be a fixed vertex in $H$. Then,
\begin{equation}
	\sum_{\substack{S \ni x ,|S|=k \\ G_S \text{ connected }}} T_k(G_S) \leq t_k^{\Delta},
\end{equation}

where, $G_S$ is the graph induced on $S$ by $G$ and, 
\begin{equation}
t_k^{\Delta} = \Delta\frac{ [(\Delta-1)(k+1)]!}{k![(\Delta-2)k + \Delta]!}.
\end{equation}

\end{theorem}

\begin{theorem}[A.Sokal \cite{SokalZeroes}]\label{Qbound} Let $Q$ be the smallest number such that,

\begin{equation}
    \inf_{\alpha > 0} \frac{1}{\alpha} \sum_{n=2}^{\infty}e^{\alpha n}Q^{-(n-1)}t_n^{(\Delta)} \leq 1.
\end{equation}

Then the choice $\alpha = 2/5$ and $Q = K\Delta = 7.963907\Delta$ satisfies the above inequality. Hence it follows that $Q \leq K\Delta = 7.963907 \Delta$.

\end{theorem}

Now we are ready to complete the proof of Theorem \ref{log-bound}. 
\begin{proof} Let $G = (V,E)$ be a graph of maximum degree $\Delta$.
Let $\mathcal{G} = (\mathcal{V}, \mathcal{E})$ denote the graph whose set of vertices is given by the set of connected subsets $S$ of $V$ such that $|S| \geq 2$ and there is an edge between $S_1$ and $S_2$ if $S_1 \cap S_2 \neq \emptyset.$ Let $X_i$ denote the set of connected subsets of $V$ of size $i$. Now we apply the above theorem for $X = \mathcal{V} = \bigsqcup_{i = 2}^{|V|} X_i$ and relation $x \sim y$ denoting that $x,y$ are disjoint in $\mathcal{G}$. \\

The generalized chromatic polynomial can be written as follows:
\begin{equation}
    \begin{split}
         & P_G(z_1, \ldots , z_q) =  \sum_{\substack{\mathcal{W} \subseteq \mathcal{V} \\ \mathcal{W} \text{ independent}}} \prod_{S_i \in S} w(S_i) \\
         & \text{where } w(S) = (z_1^{|S|} + \ldots + z_q^{|S|}) \sum_{\substack{\gamma \subseteq E,\\ (S,\gamma) \text{ connected}}} (-1)^{|\gamma|}.
    \end{split}
\end{equation}

Now we will imitate the proof of Theorem 5.1 of \cite{SokalZeroes}. We apply Theorem \ref{Dobrushin} with the choices,
\begin{equation}
    R_S = |w(S)|,
\end{equation}
and,
\begin{equation}
    A_k = \max_{x \in V}\sum_{\substack{S \ni x ,|S|=k \\ G_S \text{ connected }}} |w(S)|.
\end{equation}

This choice of $A_n$ implies that condition 1 of Theorem \ref{Dobrushin} is satisfied. 

By the definition of $\Omega$ we have,
\begin{equation}
    \begin{split}
         &|z_1^{|S|} + \ldots + z_q^{|S|}| \leq \left(2\sqrt{\frac{\Delta}{q}}\right)^{|S|-1}.
    \end{split}
\end{equation}

Thus,
\begin{equation}
       R_S =  |w(S)| \leq \left(2\sqrt{\frac{\Delta}{q}}\right)^{|S|-1}\sum_{\substack{\gamma \subseteq E,\\ (S,\gamma) \text{ connected}}} (-1)^{|\gamma|} \leq  \left(2\sqrt{\frac{\Delta}{q}}\right)^{|S|-1} T_{k}(G_S), 
\end{equation}
where $k = |S|$ and $G_S$ is the graph induced by $G$ on $S$, and $T_K(G_S)$ denotes the number of spanning trees of $G_S$. This last inequality follows from Theorem \ref{Penrose}. Thus, 

\begin{equation} 
	\begin{split}
		A_k & = \max_{x \in V}\sum_{\substack{S \ni x ,|S|=k \\ S \text{ connected }}} |w(S)| \leq  \max_{x \in V}\sum_{\substack{S \ni x ,|S|=k \\ S \text{ connected }}} \left(2\sqrt{\frac{\Delta}{q}}\right)^{|S|-1} T_{k}(G_S) \\
		& \leq \left(2\sqrt{\frac{\Delta}{q}}\right)^{k-1} \max_{x \in V}\sum_{\substack{S \ni x ,|S|=k \\ G_S \text{ connected }}}  T_{k}(G_S) \leq \left(2\sqrt{\frac{\Delta}{q}}\right)^{k-1} t_k^{\Delta}.
	\end{split}
\end{equation}
The last inequality follows from Theorem \ref{spanningbound}. Thus,
\begin{equation}
	\sum_{k=1}^{\infty} e^{\alpha k}A_k \leq \sum_{k=1}^{\infty} e^{\alpha k} \left(2\sqrt{\frac{\Delta}{q}}\right)^{k-1} t_k^{\Delta} 
\end{equation}

Let $Q$ be the smallest number such that,
\begin{equation}
    \inf_{\alpha > 0} \frac{1}{\alpha} \sum_{k=2}^{\infty}e^{\alpha n}Q^{-(k-1)}t_k^{\Delta} \leq 1.
\end{equation}

Then choosing $q$ such that, 
\begin{equation}
	\left(\sqrt{\frac{q}{4\Delta}}\right) \geq Q,
\end{equation}

gives us that, 
\begin{equation}
	\sum_{k=1}^{\infty} e^{\alpha k}A_k \leq \alpha.
\end{equation}

This gives us condition 2 of Theorem \ref{Dobrushin}, thus proving that $P_G \neq 0$ when $\left(\sqrt{\frac{q}{4\Delta}}\right) \geq Q$. By Theorem \ref{Qbound} we have $Q \leq K\Delta$. Thus, $P_G \neq 0$ when $q > 4K^2\Delta^3$.

Further, by corollary \ref{alphabound} (with $F$ being the set of edges in $G$) and Theorem \ref{Qbound} (choosing $\alpha = 2/5$) we also have that,
\begin{equation}
    |\log P_G(p_1, \ldots, p_q)| \leq 2|F|\alpha = 4|E|/5.
\end{equation}
	 
\end{proof}

\subsubsection{Taylor expansion and co-efficient bounds}

In this section we prove the bounds on the coefficients on the Taylor expansion as station in Theorem \ref{CoefBounds}. 
Using inclusion-exclusion we obtained equation (\ref{inclusion-exclusion}) for $P_G$. The following combinatorial identity is used to rewrite
the equation. Let $S_1, \ldots, S_N$ be connected subsets of $V$ and let $F(X,Y)=0$ if $X,Y$ are disjoint and -1 otherwise. Then,
\begin{equation}
        \sum_{H \in G_N} \prod_{<ij> \in H} F(S_i,S_j)  =
            \begin{cases}
                &  0 \text{ if $S_1, \ldots, S_N$ are disjoint,} \\
                & 1 \text{ otherwise, }
    \end{cases}
\end{equation}
where $G_N$ is the set of all graphs on $N$ vertices. To see this, note that the sum can be interpreted at $(1 - 1)^k$ where $k$ is the number
of pairs $(S_i, S_j)$ that are not disjoint. An intersecting pair $(S_i, S_j)$ contributes 1 to the product if $<ij>$ is not an edge in $H$
else it contributes -1. A disjoint pair $(S_i, S_j)$ contributes 0 to the sum. This gives the above identity.

Thus equation \eqref{inclusion-exclusion} can be re-written as,
\begin{equation}
          P_G(z_1, \ldots , z_q) = \sum_{N=0}^{\infty} \frac{1}{N!} \sum_{S_1, \ldots, S_N \in \mathcal{V}} \prod_{i}^{N}w(S_i) \sum_{H \in G_N} \prod_{<ij> \in H} F(S_i,S_j).
\end{equation}
The term when $N = 0$ is defined to be 1.

Using the exponential formula (\cite{EC2})one gets the Mayer expansion,
\begin{equation}\label{Taylor_exp}
        \log P_G(z_1, \ldots, z_q) = \sum_{N=1}^{\infty} \frac{1}{N!} \sum_{S_1, \ldots, S_N\in \mathcal{V}} \prod_{i}^{N}w(S_i) \sum_{H \in C_N} \prod_{<ij> \in H} F(S_i,S_j).
\end{equation}
Here $C_N$ is the set of all connected graphs on $N$ vertices.

 Let $\nu_{\mathbf{z}}(m) = z_1^m + \ldots + z_q^m$ for $2 \leq m \leq n$. The Mayer expansion is a power series of $w(S_i)$, and hence also of $\nu_{\mathbf{z}}(m)$ and the coefficients
 are independent of $q$. Theorem \ref{log-bound} tells us that $|\log P_G| \leq 4|E|/5$ on the polydisc defined by $|\nu_{\mathbf{z}}(m)| \leq \left( 2\sqrt{\frac{\Delta}{q}}\right)^{m-1}$ whenever $q > K^2\Delta^3$. This
 implies the convergence of the Mayer expansion of $P_G$ in this region. Using this we prove the bounds on the coefficients of $P_G$ as stated in Theorem \ref{CoefBounds}.

\begin{proof}
Let, $\textbf{r} = (r_1, \ldots , r_s)$ be a vector with $r_i \geq 0$. Define,
\begin{equation}
    \mathcal{M}_r(f(\nu_{\mathbf{z}}(\alpha_1), \ldots, \nu_{\mathbf{z}}(\alpha_s))) = \frac{1}{(2\pi)^s}\int_{\{(\theta_1, \ldots, \theta_s) : 0 \leq |\theta_i| \leq 2\pi \}}f(r_1 e^{i \theta_1},\ldots,r_s e^{i \theta_s} ) d\theta_1\ldots d\theta_s.
\end{equation}

Note that,
\begin{equation}
    \mathcal{M}_r\left(\frac{\nu_{\mathbf{z}}(\beta_1)\ldots \nu_{\mathbf{z}}(\beta_t)}{\nu_{\mathbf{z}}(\alpha_1)\ldots \nu_{\mathbf{z}}(\alpha_s)}\right) = 0, \text{ for } \beta \neq \alpha.
\end{equation}

Hence,

\begin{equation}
    \mathcal{M}_r\left(\frac{\log P_{G}}{\nu_{\mathbf{z}}(\alpha_1) \ldots \nu_{\mathbf{z}}(\alpha_s)}\right) =  C_{\alpha}.
\end{equation}

Also,
\begin{equation}
    \mathcal{M}_r\left(\frac{\log P_{G}}{\prod_{1}^{s} \nu_{\mathbf{z}}(\alpha_i)}\right) \leq  \frac{4|E|}{5|\prod_{1}^{s} \nu_{\mathbf{z}}(\alpha_i)|} \leq \frac{4|E|}{5\prod_{i\leq s}r_i}\\
\end{equation}

By Theorem \ref{log-bound} we know that $\log{P_G}$ converges when $q \geq K^2\Delta^3$ and $\nu_{\mathbf{z}}(\alpha_i)\leq \left(2\sqrt{\frac{\Delta}{q}}\right)^{\alpha_i -1}$. Thus, the above inequality holds when $q = K^2\Delta^3$ and $\nu_{\mathbf{z}}(\alpha_i)\leq \left(\frac{2}{K\Delta} \right)^{\alpha_i -1}$. 
Hence, using $r_i = \left(\frac{2}{K\Delta} \right)^{\alpha_i -1}$ we get,
\begin{equation}
    C_{\alpha} \leq \frac{4|E|}{5\prod_{i\leq s}\left(\frac{2}{K\Delta} \right)^{\alpha_i -1}}.
\end{equation}

Thus,

\begin{equation}
    C_{\alpha} \leq \frac{4}{5}|E|\left(\frac{K\Delta}{2} \right)^{(M-s)}, \text{ for } \alpha = (\alpha_1 \leq \ldots \leq \alpha_s) \text{ a partition of} M.
\end{equation}

Finally, we know $C_{(2,2)}$ and $C_{(3)}$ using the Mayer expansion.. Suppose vertex $i$ in $G$ has degree $d_i$. Then it can be checked from equation \ref{Taylor_exp} that,
\begin{equation}\label{coefs}
C_{(2,2)} = -\sum_{i} \binom{d_i}{2}  \text{ and } C_{(3)} \leq \sum_{i} \binom{d_i}{2}.
\end{equation}

\end{proof}

\section{Acknowledgement}

The author would like thank Prof. Persi Diaconis for introducing her to this problem and for his invaluable guidance. The author would also like to thank P\'eter Csikv\'ari for many helpful discussions. 

\bibliography{SukhadaThesis.bib}{}
\end{document}